\documentclass[leqno]{article}
\usepackage{bcp04e}

\usepackage{natbib}
\usepackage{amsmath,amssymb}

\usepackage{color}
\usepackage[dvipdf]{hyperref}

\def\E{{\mathbb E}}

\def\bP{{\bf P}}
\def\bE{{\bf E}}

\def\cG{{\mathcal G}}
\def\cF{{\mathcal F}}
\def\cH{{\mathcal H}}

\def\cG{{\mathcal G}}
\def\cQ{{\mathcal Q}}

\def\cq{\hat{q}}
\def\cp{\hat{p}}
\def\ccQ{\hat{\cQ}}

\newcounter{thm}

\newtheorem{thm}{\textbf{Theorem}}

\newtheorem{Definition}{\textbf{Definition}}

\date{\today}

\newdefin{dfn}{Definition}

\begin{document}
\mathclass{Primary 60C40; Secondary 90A46.}

\abbrevauthors{D. Ramsey and K. Szajowski}
\abbrevtitle{Correlated Equilibria in Staff Selection Problem}

\title{Correlated Equilibria in Competitive\\
        Staff Selection Problem
        }

\author{David M. Ramsey}
\address{Instytut Matematyki i Informatyki, Politechniki Wroc{\l}awskiej,\\
 Wyb\-rze\-\.ze Wys\-pia\'n\-skie\-go~27,  50-370~Wroc\-\l{}aw,
Poland\\ E-mail: ramsey@im.pwr.wroc.pl}

\author{Krzysztof Szajowski}
\address{Instytut Matematyki i Informatyki, Politechniki Wroc{\l}awskiej,\\
 Wyb\-rze\-\.ze Wys\-pia\'n\-skie\-go~27,  50-370~Wroc\-\l{}aw, Poland\\
E-mail: K.Szajowski@im.pwr.wroc.pl}

\maketitlebcp

\begin{abstract}
  This paper deals with an extension of the concept of correlated strategies
  \vfootnote{}{The idea of this paper was presented at Game Theory and Mathematical Economics, International
              Conference in Memory of Jerzy {\L}o\'s (1920 - 1998), Warsaw, September 2004
               \cite{ramsza:corrA04,ramsza:corrB04}}{}%
  to
  Markov stopping games. The Nash equilibrium approach to solving nonzero-sum
  stopping games may give multiple solutions.
  An arbitrator can suggest to each player the decision 					
  to be applied at each stage based on a joint distribution over the players' decisions.	
  This is a form of equilibrium selection. Examples of correlated equilibria in nonzero-sum
  games related to  the staff selection competition in the case of two departments are given.
  Utilitarian, egalitarian, republican and libertarian concepts of correlated equilibria selection are used.
\end{abstract}

\section{Introduction}
In this paper an alternative approach to the staff selection competition
in the case of two departments considered by Baston and Garnaev \cite{basgar04:staff} is proposed.
The formulation of the problem in Baston and Garnaev \cite{basgar04:staff} is as follows.
Two departments in an organisation are each seeking to make an appointment
within the same area of expertise. The heads of the two departments together
interview the applicants in turn and make their decisions on one applicant
before interviewing any others. If a candidate is rejected by both departmental heads,
the candidate cannot be considered for either post at a later date. When both heads decide to make
an offer, they consider the following possibilities.
\begin{enumerate}
\item\label{casea} The departments are equally attractive, so that an applicant has no preference between them;
\item\label{caseb} One department can offer better prospects to applicants, who will always choose that
           department.
\end{enumerate}
The departmental heads know that there are precisely $N$ applicants and that each applicant
has a level of expertise which is random. It is assumed that the interview process
enables the directors to observe these levels of expertise, which form a sequence of
\emph{i.i.d} random variables from a continuous distribution. If no
appointment is made to a department from these $N$ applicants, then the department
will suffer from a shortfall of expertise. Game \ref{caseb} has one Nash equilibrium,
 which can be used as the solution to the problem. Game \ref{casea} has many Nash equilibria.
This raises the question of equilibrium selection. Baston and Garnaev
\cite{basgar04:staff}
interpreted such a variety of Nash equilibria solutions as a way of modelling different dynamics
within the organisation, which can result in various outcomes during the conscription process.
If one departmental head is aggressive and one passive, we might expect a different outcome
to the one in which both are of a similar temperament. When both have a
similar temperament one expects a symmetric
strategy and value, but when they have different temperaments
one should expect an asymmetric equilibrium and value. The different character
of heads is modelled by the notion of a Stackleberg leader.
Also, the difference in the level of complication of equilibria might
 also be an argument justifying this approach to equilibrium
selection. It is shown that these non-symmetric equilibria have the advantage that the players use pure strategies,
whereas at the symmetric equilibrium, the players are called upon to employ specific actions
with complicated probabilities.

The  staff selection problem presented above is closely related to the best choice problem (BCP).
There are some potential real applications of decision theory which strengthen the motivation of the BCP
(the one decision maker problem). One group of such problems are models of many important business decisions, such
as choosing a venture partner, adopting technological innovation, or hiring an employee using a sequential decision
framework (see Stein, Seale and Rapoport~\cite{stesearap03:BCP},
Chun~\cite{chupla92:selection,chu96:weighted,chu99:group,chu00:search}).
Others are an experimental investigations of the ,,secretary problem'', which compare the optimal policy from
the mathematical model with behaviour of human beings
 (see Seal and Rapoport~\cite{searap97:exper,searap00:unknown}).
We have not found any such investigation for BCP games. It could be that the
theoretical results are not complete enough to start
applied and experimental research.

In spite of the long history of BCP and its generalisations presented in review papers by
Freeman \cite{fre83:review}, Ferguson \cite{fer89:who}, Rose \cite{rose82:survey},
 Samuels \cite{sam91:secpro}, there are also competitive versions, on which researchers' attention has been focused
(see Sakaguchi \cite{sak89:multi,sak95:review} for review papers). Let us briefly
recall the main game theoretic models of BCP.
Enns and Ferenstein \cite{ennfer85:horse}, Enns, Ferenstein and Sheahan
\cite{ennfershe86:curious} solved
a non-zero sum game related to BCP. Some important mathematical results
 related to the problem, posed in this paper,
were proven many years later by Bruss and Louchard \cite{brulou98:horse}. The full information version of the game
was solved by Chen, Rosenberg and Shepp \cite{cherosbur97:two}.
The relation between players is as follows.
The players have numbers: $1$ and $2$.
When an item appears then Player $1$ always has the first opportunity to decide
whether to hire the applicant or not (unless
she has hired one already). One can say that Player $1$ has priority. If
Player $1$ does not hire the current applicant, then Player $2$ can decide
whether to hire the applicant
or not (unless she has hired one already).
 If neither player hires the current applicant,
they interview the next applicant. The interview process continues until both players
have hired an applicant.
A hired applicant does not hesitate and accepts an offer without any delay or
 additional conditions. The games in this group of papers
have the same strategic scheme as in Game \ref{caseb}.

The concept of equal priority of the players in the selection process
in a model of a non-zero-sum game related to BCP
was introduced by Fushimi \cite{fus81:competitive}. Szajowski \cite{sza94:markov}
 extended this model
to permit random priority. Ramsey and Szajowski \cite{ramsza01:random,ramsza05:bilateral}
 considered a mathematical model of
competitive selection with random priority and random acceptance of the offer (uncertain employment)
by candidates. Uncertain employment is a source of additional problems, which are solved as follows.
At each moment $n$ the candidate is presented to both players. If neither player has yet obtained an object then:
\begin{description}
\item [(i)] if only one of them would like to accept the state,
            then he tries to take it. In this case the random
            mechanism assigns the availability of the state (which
            can depend on the player and the moment of decision $n$);
\item [(ii)] if both of them are interested in this state, then
            the random device chooses the player who will first solicit
            the state. The availability of the state is the same as in the
            situation when only one player wants to take it. If the chosen
            player obtains the state, he stops searching;
\item [(iii)] If this state is not available to the player chosen by
            the random device, then the observed state at moment $n$ is lost to both players.
            Both players continue searching by inspecting the next state.
\end{description}

When one player has obtained a candidate the other player continues searching alone. If this player wishes to
accept a candidate, the probability that it is available to him is the same as in point (i) above.

When a non-zero-sum game does not have
a unique Nash equilibrium, then
 communication between the players
would be useful in deciding which equilibrium should be played. Using the idea of correlated strategies introduced by
\cite{aum74:subject},
the set of possible strategies is extended to the set of correlated stopping times and
the actions undertaken by the players are correlated.

Little research has been carried out on the role of communication between players in stopping games. \cite{sol01:corr}
and \cite{sol02:corrstoch} consider correlated equilibria in general dynamic games. The form of correlation
is not unique. The approach applied here is based on a generalisation of randomised
 stopping times.
Various additional criteria used by the players to correlate their actions restrict the set of possible solutions.
These criteria are based on those used in \cite{grehal03a:learning}, which
resemble ideas of solutions of cooperative games presented in \cite{tho94:bargain}.

Strategies of staff selection based on the construction of correlated
strategies according to various selection
criteria are presented in the setting adopted by Baston and Garnaev
\cite{basgar04:staff}.
Correlated strategy selection was proposed by the authors in
\cite{ramsza03:corr}.

The construction of correlated equilibria in stopping games is based on the
concept of correlated equilibria in two-by-two
bimatrix games. The geometry of correlated equilibria in bimatrix games is
 described by Calv\'o-Armengol \cite{cal04:corre}.

\section{\label{corrstop1}Correlated equilibria in stopping games }
\cite{aum74:subject} introduced a correlation scheme in randomised
strategies for non-zero-sum games extending the concept of Nash equilibrium. Using
this approach some process of preplay communication is needed to realise such a strategy.
Aumann's approach has been extended in various manners (eg see
\cite{for86,germou78:corr,moulin,now93:corr,tolhaulei86:stoch}). The process of
adapting correlated equilibria to stopping games starts from the idea of
correlated stopping times.
\begin{Definition}\label{corrstrat1}
A random sequence $\cq =\{(q^1_n,q^2_n,q^3_n)\}$ such that, for each $n$,
\begin{description}
\item[(i)]   $q^i_n$ is adapted to $\cF_n$ for $i=1,2,3$;
\item[(ii)]  $0\leq q^1_n\leq q^2_n\leq q^3_n \leq 1$ a.s.
\end{description}
is called a correlated stopping strategy. The set of all such sequences will be
denoted by $\ccQ^N$.
\end{Definition}

Let $A_1,A_2,\ldots,A_N$ be a sequence of i.i.d. r.v. with uniform distribution
on $[0,1]$ and independent of the Markov process $(X_n,\cF_n,\bP_x)_{n=0}^N$.
Denote $\vec{q}_n=(q^1_n,q^2_n,q^3_n)$. Correlated stopping times are
pairs $(\lambda^1(\cq), \lambda^2(\cq))$ of Markov stopping times with respect to the $\sigma$-fields
$\cH_n=\sigma\{\cF_n,\ \{A_1,A_2,...,A_n \}\}$ defined by the strategy
$\cq = (\vec{q}_n)\in \ccQ^N$ as follows:
\begin{equation}\label{corrstopt1}
\lambda^1(\cq) = \inf\{0\leq n\leq N: A_n\leq q^2_n\}
\end{equation}
and
\begin{equation}\label{corrstopt2}
\lambda^2(\cq) = \inf\{0\leq n\leq N: A_n\leq q^1_n \mbox{ or }
                    q^2_n< A_n\leq q^3_n\}.
\end{equation}
The strategy $\cq$ will be called the correlation profile and it defines the pair
of stopping times $(\lambda^1(\cq),\lambda^2(\cq))$.

	In intuitive terms, the vector $\vec{q}_n=(q^1_n ,q^2_n ,q^3_{n})$ defines the joint distribution of the actions
 	taken by the players at moment $n$:
 	with probability $q^1_n$ both players choose the action "stop", with probability $q^2_n-q^1_{n}$ Player 1 stops
 	and Player 2 chooses the action "continue", with probability $q^3_n-q^2_n$ Player 1
   	continues and Player 2 stops and with probability $1-q^3_n$ both players continue. A correlated strategy $\hat{q}$
    	is assumed to be defined by preplay communication  between the players (either before the start of the game or before each
	decision) with the possible aid of an "external judge". If communication only takes
	place before the game commences, then such a correlation is said to be a stationary correlation device.
	If communication may occur at each decision point, then such a correlation is said to be an extensive (autonomous)
	correlation device (see \cite{sol02:corrstoch}). In general, we consider extensive correlation
	devices. The form of the correlated strategy is known to both players.

If one player carries out the actions suggested by the external judge with the aid of the
appropriate lottery and
the other player departs from the suggested action a formal construction of the possible strategies and the
calculation of the expected gains should be done.

Let $\hat{p}=(p_{1},p_{2},\ldots,p_{N})$ be a sequence in the
unit interval. If Player $i$ departs from the correlation profile $\cq$, then the strategy of the other player is based
on the marginal correlated profile $\cq_{-i}$ and the strategy of Player $i$ is defined by $\hat{p}_i=\hat{p}$.
Denote $\tau^i((\cp_i,\cq_{-i}))=\tau^i(\cp_i)=\inf\{0\leq n\leq N: A^{'}_n\leq p_{n}\}$,
where $(A^{'}_n)_{n=1}^N$ is a sequence of i.i.d. r.v. with uniform distribution
on $[0,1]$, independent of $(A_n)_{n=1}^N$ and independent of the Markov process $(X_n,\cF_n,\bP_x)_{n=0}^N$.
Denote $\bar{G}_i(\cq)=G_i(\lambda^1(\cq)\wedge \lambda^2(\cq),X_{\lambda^1(\cq)\wedge \lambda^2(\cq)})$ and
$\bar{G}_i((\cp_i,\cq_{-i}))=
G_i(\tau^i(\cp_i)\wedge \lambda^{-i}(\cq_{-i}),X_{\tau^i(\cp_i)\wedge \lambda^{-i}(\cq_{-i})})$.
The expected payoffs are defined as $\hat{G}_i(x,\cq)=\bE_x\bar{G}_i(\cq)$ and $\hat{G}_i(x,(\cp_i,\cq_{-i}))=
\bE_x\bar{G}_i((\cp_i,\cq_{-i}))$, respectively.

\begin{Definition}\label{defCE}
A correlated stopping strategy $\cq^*\in\ccQ^N$ is called a correlated equilibrium point of $\cG_m$, if
\begin{equation}\label{defCE1}
\hat{G}_i(x,\cq^{*})\geq \hat{G}_i(x,(\cp_i,\cq^{*}_{-i}))
\end{equation}
for every $x\in\E$, $\hat{p}$ and $i=1,2$.
\end{Definition}

This is a definition of a correlated equilibrium in the normal form of the game.
It should be noted that a stronger notion of correlated equilibrium can be introduced
by requiring that the correlation must define an equilibrium in each restricted game
where $n$ steps remain $(1\leq n\leq N-1)$.

\section{\label{corrselect}Selection of a Correlated Equilibrium}%
Since the set of Nash equilibria is a subset of the set of correlated equilibria, it is clear that whenever the
problem of the selection of a Nash equilibria exists, the problem of the selection of a correlated equilibrium
also exists. However, the notion of correlated equilibrium assumes that communication takes place.
Such communication can be used to define the criteria used by players to select a correlated equilibrium.
 We now formulate various criteria for selecting a correlated equilibria. These criteria select
subsets of $\mathbb{CE}$. The concepts which are used here do not come from the concepts of solution to Nash's
 problem of cooperative bargaining. These concepts were used by Greenwald and Hall \cite{grehal03a:learning}
  for computer learning of equilibria in Markov games.

\begin{Definition}\label{corrselect1}
Let us formulate four different selection criteria for correlated equilibria in a stopping game.
\begin{enumerate}
\item\label{util} A utilitarian correlated equilibrium
is an equilibrium constructed recursively in such a way that at each stage
$n=N-1,N-2,\ldots ,1$
the sum of the values of the game to the players is maximised
given the equilibrium calculated for stages $n+1,n+2,\ldots ,N$ is played..

\item\label{egal} An egalitarian correlated equilibrium
is an equilibrium constructed recursively in such a way that at each
stage $n=N-1,N-2,\ldots ,1$ the minimum value is maximised
given the equilibrium calculated for stages $n+1,n+2,\ldots ,N$ is played.

\item\label{repub} A republican correlated equilibrium
is an equilibrium constructed recursively in such a way that at each stage
$n=N-1,N-2,\ldots ,1$ the maximum value is maximised
given the equilibrium calculated for stages $n+1,n+2,\ldots ,N$ is played.

\item\label{libert} A libertarian $i$ correlated equilibrium
is an equilibrium constructed recursively in such a way that at each stage
$n=N-1,N-2,\ldots ,1$
 the value of the game to Player $i$ is maximised
given the equilibrium calculated for stages $n+1,n+2,\ldots ,N$ is played.
\end{enumerate}
\end{Definition}

\begin{thm}\label{thcorreq1a}
The set of correlated equilibrium points satisfying any one of the given criteria above is not empty.
\end{thm}

\section{\label{staffselect1}One and two applicant games with no candidate preferences}%
Let us assume that the cost of not selecting an applicant is $c$.
This is the cost of a shortfall of
expertise in a department. If a director selects an applicant with expertise
$\xi_i=x$, the department gains $x$. Let
us assume that the candidates have \emph{i.i.d.} expertise $\xi_i$ with
 uniform distribution on $[0,1]$. If there is
only one candidate, then the selection process will end with value
$d=\frac{1}{2}\bE\xi_1-\frac{1}{2}c=\frac{1-2c}{4}$
to both players (both want to select and the probability of winning is $\frac{1}{2}$ for both of them).
Let $b = \max \{ 0, \frac{1-2c}{4} \}$.
\subsection{Correlated equilibria of the two stage game}
When there are two candidates, then we have a two stage game. The subgame perfect
Nash equilibria at the stage when
the first candidate is interviewed will be considered. The payoff bimatrix
$M_2(x)$ is of the form (see \cite{basgar04:staff}):
\begin{equation}
M_2(x)=\begin{array}{cc}
                    &\mbox{\hfil s\hspace{5em} f\hfil}\\
        \begin{array}{c}
        s\\
        f
        \end{array}&\multicolumn{1}{c}{
                    \left(\begin{array}{cc}
                	(\frac{(x+\frac{1}{2})}{2},\frac{(x+\frac{1}{2})}{2})&(x,\frac{1}{2})\\
                	(\frac{1}{2},x)&				      (d,d)
                    \end{array}
                    \right)}
       \end{array}
\end{equation}
The game has one pure Nash equlilibrium, $(s,s)$, for $x\ge \frac{1}{2}$
 and for $x\le b$ has one pure Nash
equilibrium $(f,f)$. However, for $x\in[b,\frac{1}{2}]$ there are two
asymmetric pure Nash equilibria and one
symmetric Nash equilibrium in mixed strategies. Without extra assumptions it is not clear
which equilibrium should be played.
Baston and Garnaev~\cite{basgar04:staff} have proposed that if the
players have a similar character, then the symmetric solution
should be played. In the non-symmetric case the idea
of Stackleberg equilibrium can be adopted.
It is assumed that the first player will be the Stackleberg leader and the $1$-Stackleberg equilibrium
is the solution of the problem selected.

We will use an extensive communication device to construct correlated equilibria.
 In general,
 correlated
equilibria are not unique. Usually the set of correlated equilibria contain the
convex hull of Nash equilibria. However,
natural selection criteria can be proposed
and the possibility of preplay communication and use of an arbitrator solve
the problem of solution selection. The players just specify the criterion.
Such criteria are formulated in
Section \ref{corrselect}. The set of
solutions which fulfil one of the points \ref{util}-\ref{libert} in definition \ref{corrselect1} are not empty.

For $M_2(x)$, when $x\in [b,\frac{1}{2}]$ the set of correlated equilibria is a polytope with five vertices.
Let $\alpha=\frac{1}{2}\frac{x-\frac{1}{2}}{d-x}$ and $\gamma=2\frac{d-x}{x-\frac{1}{2}}$ and let us denote
$\mu=(\mu_{ss},\mu_{ff},\mu_{fs},\mu_{sf})$. The polytope of correlated equilibria for the considered game has
the five vertices given in Table \ref{opvertices} (see Peeters and Potters \cite{peepot99:structure}).
\begin{table}
\[
\begin{array}{|l||c|c|c|c|}\hline\hline
\mu				&\mu_{ss}	&\mu_{ff}	&\mu_{fs}	&\mu_{sf}   	\\ \hline\hline
\mu_C^{\ast}(\alpha,\gamma)     &0		&0		&1		&0		\\ \hline
\mu_D^{\ast}(\alpha,\gamma)     &0		&0		&0		&1		\\ \hline
\mu_E^{\ast}(\alpha,\gamma)     &\frac{\gamma}{1+\gamma+\alpha\gamma}&0&
                                 \frac{1}{1+\gamma+\alpha\gamma}&\frac{\alpha\gamma}{1+\gamma+\alpha\gamma}\\ \hline
\mu_F^{\ast}(\alpha,\gamma)     &0				&\frac{\alpha}{1+\alpha+\alpha\gamma}&
				 \frac{1}{1+\alpha+\alpha\gamma}&\frac{\alpha\gamma}{1+\alpha+\alpha\gamma}\\ \hline
\mu_G^{\ast}(\alpha,\gamma)     &\frac{\gamma}{(1+\alpha)(1+\gamma)}&\frac{\alpha}{(1+\alpha)(1+\gamma)}&
 				 \frac{1}{(1+\alpha)(1+\gamma)}&\frac{\alpha\gamma}{(1+\alpha)(1+\gamma)}  \\ \hline\hline
\end{array}
\]
\caption{\label{opvertices}The five vertices of the correlated equilibrium polytope.}
\end{table}
The value at each vertex will be calculated.
\begin{description}
\item[(C)]  The values of the game to the players at vertex $C$ are denoted
by $v_1^{(C)}$ and $v_2^{(C)}$.
\begin{eqnarray}
\label{v1C}%
v_1^{(C)}&=&\int_0^bbd\!x+\int_b^\frac{1}{2}\frac{1}{2}d\!x+\frac{1}{2}\int_\frac{1}{2}^{1}(x+\frac{1}{2})d\!x
          =b^2-\frac{1}{2}b+\frac{9}{16}\\
\label{v2C}%
v_2^{(C)}&=&\int_0^bbd\!x+\int_b^\frac{1}{2}xd\!x+\frac{1}{2}\int_\frac{1}{2}^{1}(x+\frac{1}{2})d\!x
          =\frac{1}{2}b^2+\frac{7}{16}
\end{eqnarray}
When Player $1$ takes the role of Stackleberg leader his expected gain is $v_1^{(C)}$,
 while the Stackleberg follower has
$v_2^{(C)}$ (see \cite{basgar04:staff}).

\item[(D)] The values at vertex $D$ can be obtained from those at vertex $C$,
because matrix $M_2(x)$ is symmetric.
\begin{eqnarray}
\label{v1D}%
v_1^{(D)}&=&\frac{1}{2}b^2+\frac{7}{16}\\
\label{v2D}%
v_2^{(D)}&=&b^2-\frac{1}{2}b+\frac{9}{16}
\end{eqnarray}

\item[(E)] The expected gain of the players at correlated equilibrium $E$ given the
expertise of the candidate
$x\in [b,\frac{1}{2}]$ is of the form.
\begin{eqnarray}
\label{w1E}%
w_1^{(E)}&=&(x+\frac{1}{2})\frac{x-\frac{1}{2}}{2(d-\frac{1}{2})}+\frac{1}{2}(x+\frac{1}{2})\frac{d-x}{d-\frac{1}{2}}\\
\nonumber&=&\frac{1}{2}(x+\frac{1}{2})\\
\label{w2E}%
w_2^{(E)}&=&\frac{1}{2}(x+\frac{1}{2}).
\end{eqnarray}
The value of the two-stage game to the players at vertex $E$ is
\begin{eqnarray}
\label{v1E}%
v_1^{(E)}&=&v_2^{(E)}=\int_0^bbd\!x+\frac{1}{2}\int_b^1(x+\frac{1}{2})d\!x=\frac{3}{4}b^2-\frac{1}{4}b+\frac{1}{2}.
\end{eqnarray}
The values at these three vertices are such that $v_1^{(D)}< v_1^{(E)}< v_1^{(C)}$.

\item[(F)] This correlated equilibrium is of the form: $\mu_{ss}=0$ and
\begin{eqnarray*}
\mu_{ff}&=&\frac{x-\frac{1}{2}}{4d-3x-\frac{1}{2}}\\
\mu_{sf}&=&\frac{2(d-x)}{4d-3x-\frac{1}{2}}\\
\mu_{fs}&=&\mu_{sf}.
\end{eqnarray*}
The expected gain of the players at correlated equilibrium $F$
 given the expertise of candidate
$x\in [b,\frac{1}{2}]$ is
\begin{eqnarray}
\label{w1F}%
w_1^{(F)}&=&w_2^{(F)}=\frac{d(x-\frac{1}{2})+2(d-x)(x+\frac{1}{2})}{4(d-x)+x-\frac{1}{2}}\\
\nonumber&=&\frac{1}{2}(x+\frac{1}{2})+\frac{(x-\frac{1}{2})(d-\frac{x}{2}-\frac{1}{4})}{4d-3x-\frac{1}{2}}\\
\nonumber&\le&\frac{1}{2}(x+\frac{1}{2})
\end{eqnarray}
for $x\in [b,\frac{1}{2}]$. The value of the two-stage game to the
 players at vertex $F$ is
\begin{eqnarray}
\label{v1F}%
v_1^{(F)}&=&v_2^{(F)}=v_1^{(E)}+\int_b^\frac{1}{2}\frac{(x-\frac{1}{2})(d-\frac{x}{2}-\frac{1}{4})}{4d-3x-\frac{1}{2}} d\!x\\
\nonumber & < & v_1^{(E)}.
\end{eqnarray}

\item[(G)] This correlated equilibrium (the Nash equilibrium in mixed strategies)
is of the form:
\begin{eqnarray*}
\mu_{ss}&=&\frac{4(d-x)^2}{(2d-x-\frac{1}{2})^2}\\
\mu_{ff}&=&\frac{(x-\frac{1}{2})^2}{(2d-x-\frac{1}{2})^2}\\
\mu_{sf}&=&\frac{2(d-x)(x-\frac{1}{2})}{(2d-x-\frac{1}{2})^2}\\
\mu_{fs}&=&\mu_{sf}.
\end{eqnarray*}
The expected gain of the players at correlated equilibrium $G$ given the
expertise of the candidate
$x\in [b,\frac{1}{2}]$ is
\begin{eqnarray}
\label{w1G}%
w_1^{(G)}&=&w_2^{(G)}=\frac{2(d-x)^2(x-\frac{1}{2})+2(d-x)(x+\frac{1}{2})(x-\frac{1}{2})+d(x-\frac{1}{2})^2}%
                       {(2d-x-\frac{1}{2})^2}\\
\nonumber&=&\frac{1}{2}(x+\frac{1}{2})+\frac{(x-\frac{1}{2})^2[
d-\frac{1}{2}(x+\frac{1}{2})]}{(2d-x-\frac{1}{2})^2}\\
\nonumber&\le&\frac{1}{2}(x+\frac{1}{2})
\end{eqnarray}
for $x\in [b,\frac{1}{2}]$. The value of the two-stage game to the
 players at vertex $G$ is
\begin{eqnarray}
\label{v1G}%
v_1^{(G)}&=&v_2^{(G)}=v_1^{(E)}
          +\int_b^\frac{1}{2}\frac{(x-\frac{1}{2})^2[d-\frac{1}{2}(x+\frac{1}{2})]}
{(2d-x-\frac{1}{2})^2} d\!x\\
\nonumber & < & v_1^{(E)}.
\end{eqnarray}

\end{description}

\subsection{Selection of equilibria in the two stage game}
Let us apply the selection criteria on the set of correlated equilibria of the two
stage game.  We thus define a linear programming problem, in
which the objective function is defined by the criterion
and the feasible set is
the set of vectors $\mu$ defining a correlated equilibrium. Hence to find a
solution, we compare
the appropriate values at each vertex of the correlated equilibria polytope described
in the previous section.

It should be noted that
when either the republican or
egalitarian criterion is used, the solution is given by the appropriate solution
from one of two linear programming problems. In these cases the
two linear programming problems are:

\begin{enumerate} \item[1)] Maximise $v_{1}$ given the equilibrium
constraints and the constraint $v_{1}\leq v_{2}$ when the egalitarian condition is
used or $v_{1}\geq v_{2}$ when the republican condition is used.

\item[2)] Maximise $v_{2}$ given the equilibrium
constraints and the constraint $v_{2}\leq v_{1}$ when the egalitarian condition is
used or $v_{2}\geq v_{1}$ when the republican condition is used.
\end{enumerate}

From the symmetry of the game the hyperplane $\mu_{fs}-\mu_{sf}=0$ splits the
set of correlated equilibria into the two feasible sets for these problems and
$\mu = (0,0,\frac{1}{2},\frac{1}{2})$ becomes a vertex of the feasible set in each of
the problems. We call this vertex $H$. This vertex replaces vertex $C$ or vertex $D$
depending on the additional constraint. We have
\begin{equation}
v_{1}^{(H)} = v_{2}^{(H)} = \frac{v_{1}^{(C)}+v_{1}^{(D)}}{2} = v_{1}^{(E)}
\label{v1H}
\end{equation}

\subsubsection{Libertarian equilibria}
From (\ref{v1C})--(\ref{v1G}) it follows that the maximal game value
for the first player is guaranteed at vertex
$(f,s)$ and for the second player at $(s,f)$.
 It means that $\delta_{L1}^\star=(f,s)=C$ is the libertarian $1$
and $\delta_{L2}^\star=(s,f)=D$ is the libertarian $2$ correlated equilibrium.
In relation to the solutions presented by Baston and Garnaev, the libertarian $i$
equilibrium corresponds to the Stackleberg solution at which Player $i$ takes the
role of the Stackleberg leader.

\subsubsection{Egalitarian equilibria}
Let us denote $v^\delta=\min_{i\in \{1,2\}}v_i^\delta$.
We are looking for $\delta_E^\star$ such that $v^{\delta_E^\star}=\max_\delta v^\delta$. For $\delta\in\{E,F,G,H\}$ we have
$v_1^\delta=v_2^\delta$, $v_1^{(F)}< v_1^{(E)}=v_{1}^{(H)}$ and $v_1^{(G)}< v_1^{(E)}$.
For $\delta\in\{C,D\}$ the minimal values
are $v^{(C)}=v_2^{(C)}$ and $v^{(D)}=v_1^{(D)}$. Moreover, $v_2^{(C)}=v_1^{(D)}< v_1^{(E)}$.
Therefore $E$ and $H$ define egalitarian equilibria and $v^{\delta_E^\star}=v_1^{(E)}$.
It follows that any linear combination of these equilibria $pE+(1-p)H$, where
$p\in [0,1]$ defines an
egalitarian equilibrium. It should be noted that $H$ is an intuitively pleasing
solution, since it corresponds to a solution in which the players observe the toss of a
coin and if heads appears Player 1 acts as the Stackleberg leader, otherwise
Player 2 plays this role. This is one of the solutions considered by Baston and
Garnaev. At any of the other solutions the arbitrator must send signals to each of
the players separately in order to obtain the appropriate correlation. It should be noted
that the value of the game to the players is independent of the egalitarian
solution adopted.

\subsubsection{Republican equilibria}
Let us denote $V^\delta=\max_{i\in \{1,2\}}v_i^\delta$.
Similar consideration of the vertices as made in the case of
egalitarian equilibria leads to conclusion that the republican equilibria are
$\delta_R^\star\in\{C,D\}$ and $V^{\delta_R^\star}= v_1^{(C)}=v_2^{(D)}$. These
are the only two solutions, since they are the unique solutions of the
two appropriate linear programming programmes described above and correspond
to the Stackleberg solutions.

\subsubsection{Utilitarian equilibria}
Let us denote $v_{+}^\delta=v_1^\delta+v_2^\delta$.
We have $v_{+}^{(C)}=v_{+}^{(D)}=\frac{3}{2}b^2-\frac{b}{2}+1=2v_1^{(E)}$.
 Since $2b\le x+\frac{1}{2}$, it follows that
$v_{+}^{(C)}>v_{+}^{(F)}$ and $v_{+}^{(C)}>v_{+}^{(G)}$.
Hence, $C,D$ and $E$ are utilitarian equilibria. It follows that any linear
combination $pC+qD+rE$ ($p,q,r \geq 0, p+q+r=1$) defines a utilitarian equilibrium.
$v_{+}^{\delta^{\star}_U}=v_{+}^{(C)}=v_{+}^{(D)}=v_{+}^{(E)}$. It should be noted
that $H$ is a linear comibination of these three vertices with
$p=q=\frac{1}{2},r=0$. Also, the value of the game to the players is dependent on the utilitarian
equilibrium played.

\section{Selection of equilibria in the multi-stage game}

We define correlated equilibria by recursion as a series of correlated equilibria in
the appropriately defined matrix games. The correlated strategy used when both players
are deciding whether to accept or reject the $n$-th last candidate is given by
$\mu_{n} = (\mu_{n,ss},\mu_{n,ff},\mu_{n,fs},\mu_{n,sf})$.
The game played on observing the $n$-th last candidate is given by
\[
M_n(x)=\begin{array}{cc}
                    &\mbox{\hfil s\hspace{5em} f\hfil}\\
        \begin{array}{c}
        s\\
        f
        \end{array}&\multicolumn{1}{c}{
                    \left(\begin{array}{cc}
                	(\frac{x+u_{n-1}}{2},\frac{x+u_{n-1}}{2})&(x,u_{n-1})\\
                	(u_{n-1},x)&				      (v^{\pi}_{n-1},w^{\pi}_{n-1})
                    \end{array}
                    \right)} ,
       \end{array}
\]
where $u_{n}$ is the optimal expected reward of a lone searcher with $n$ candidates
remaining (see [2]) and $v^{\pi}_{n},w^{\pi}_{n}$ are the values of
the $n$-stage game  to
Players 1 and 2, respectively, when the equilibrium $\pi$ is played. From the form of the payoff matrix it can be seen
that $(s,s)$ is the unique Nash equilibrium when $x >u_{n-1}$. Similarly, $(f,f)$ is
the unique Nash equilibrium when $x< \min \{ v^{\pi}_{n-1}, w^{\pi}_{n-1} \} $.

\subsection{Libertarian equilibria}

First we consider $N=3$. From the calculations made for $N=2$, it follows
that $v^{L1}_{2}>w^{L1}_{2}$. Considering the payoff matrix $(f,s)$ is the unique Nash
equilibrium for $v^{L1}_{2}<x<w^{L1}_{2}$ and both $(f,s)$ and $(s,f)$ are pure
Nash equilibrium for $v^{L1}_{2}<x<u_{2}$. There is also an equilibrium in mixed
strategies on this interval. Thus, we only need to consider equilibrium selection
for $v^{L1}_{2}<x<u_{2}$. Since the payoff matrix is now longer symmetric,
the vertices
of the polytope defining the set of correlated equilibrium are of a different form.
However, since $(f,s)$ is a Nash equilibrium, $\mu_{3} = (0,0,1,0)$ is a
vertex of this polytope. For $v^{L1}_{2}<x<u_{2}$, it can be seen that
$u_{2}$ is the maximal payoff in the payoff matrix. It follows that
$\mu_{3}$ is the vertex that strictly
maximises the
expected payoff of Player 1 and thus uniquely defines the libertarian 1 equilibrium. It
follows that $v^{L1}_{3}>w^{L1}_{3}$ and hence $M_{4}(x)$ is of a similar form to
$M_{3}(x)$. By iteration it follows that Player 1 plays the
role of the Stackleberg leader at the libertarian 1 solution.
Analogously, Player 2 plays the role of the Stackleberg leader at the libertarian 2
solution. For the value functions see [2].

\subsection{Egalitarian equilibrium}

It will be shown by induction that for $N\geq 3$ an egalitarian equilibrium is of
the same form as for $N=2$. Suppose
that $v^{E}_{n-1}=w^{E}_{n-1}$. The coordinates of the vertices of the polytope
describing the set of correlated equilibria is of the form given in
Table \ref{opvertices} with $\alpha = \frac{u_{n-1}-x}{2(x-v_{n-1})}$ and
$\gamma = \frac{2(x-v_{n-1})}{u_{n-1}-x}$. Considering the values of the game at
these vertices when $x\in [ v_{n-1},u_{n-1} ]$, the egalitarian
criterion is satisfied at
vertices $E$ and $H$. It follows that $v^{E}_{n}=w^{E}_{n}$ and any
linear combination of $E$ and $H$ defines an egalitarian equilibrium. Since
$v^{E}_{2}=v^{E}_{2}$ it follows by induction that an egalitarian equilibrium is of
the required form. In particular, the equilibrium obtained by deciding who plays the
role of Stackleberg leader based on the result of a coin toss defines an
egalitarian equilibrium.

\subsection{Republican equilibria}

Suppose libertarian 1 is taken to be the
republican equilibrium for the last 2 stages. For $N=3$ the calculations are similar to
the calculations made for the libertarian 1 equilibrium. It can be shown that the
libertarian 1 equilibrium again maximises the maximum value. Using an iterative
argument, it can be shown that the libertarian 1 equilibrium is a republican
equilibrium. By the symmetry of the game it follows that the libertarian 2
equilibrium is also a republican equilibrium.

\subsection{Utilitarian equilibria}

Unfortunately, the value function of a utilitarian equilibrium for $N=2$ is not
uniquely defined. In order to find a "globally optimal" utilitarian equilibrium, we
cannot use simple recursion. From the form of the payoff matrix it can be
seen that when $\max \{ v_{n-1},w_{n-1} \} < x < u_{n-1}$
the maximum sum of payoffs is $x+u_{n-1}$. This is obtained when at least one of the
players accepts the candidate. Such a payoff is attainable at a correlated equilibrium,
since $(f,s)$ and $(s,f)$ are correlated equilibrium. It follows from the
definition of a utilitarian equilibrium that $\mu_{n,ff}=0$ when
$\max \{ v_{n-1},w_{n-1} \} < x < u_{n-1}$.

\begin{thm}\label{thmutil}
The libertarian equilibria are the only globally optimal utilitarian equilibria for
$N\geq 3$ (ignoring strategies whose actions differ from those defined by one of these strategies on a set with
probability measure zero).
\end{thm}

{\bf Proof} First we show that among the set of utilitarian
equilibria the minimum value is minimised at the libertarian equilibria for $N\geq 2$.
Considering
the values of the game at the vertices of the set of utilitarian
correlated equilibria when $N=2$ (obtained by adding the additional condition that
$\mu_{2,ff}=0$ for $b<x<\frac{1}{2}$),
the minimum value is minimised at the two libertarian equilibria. From the form of
the two linear programming problems that define this minimisation problem, it follows
that these solutions are the only such solutions.

By symmetry $w^{L1}_{n}
= v^{L2}_{n}$. Set $k^{\pi}_{n}=\min \{ v^{\pi}_{n}, w^{\pi}_{n} \}$.
Assume that $w^{L1}_{n}<k^{\pi }_{n}$ for all $\pi \neq L1, L2 $. We have
\begin{eqnarray*}
w^{L1}_{n+1} & = & \int_{0}^{w^{L1}_{n}} w^{L1}_{n} dx +
\int_{w^{L1}_{n}}^{u_{n}} x dx
+ \frac{1}{2} \int_{u_{n}}^{1} (x+u_{n}) dx \\
k^{\pi}_{n+1} & \geq & \int_{0}^{k^{\pi}_{n}} k^{\pi}_{n} dx +
\int_{k^{\pi}_{n}}^{u_{n}} g^{\pi}_{n}(x) dx
+ \frac{1}{2} \int_{u_{n}}^{1} (x+u_{n}) dx ,
\end{eqnarray*}
where $g^{\pi}_{n}(x)$ is the expected reward of such a player given that
$x\in [k^{\pi}_{n},u_{n}]$. From the condition that $\mu_{n,ff}=0$ it follows
that $g^{\pi}_{n}(x)\geq x, \forall x\in [k^{\pi}_{n},u_{n}]$ and
\begin{displaymath}
k^{\pi}_{n+1} - w^{L1}_{n+1} \geq (k^{\pi}_{n})^{2}-(w^{L1}_{n})^{2} -
\int_{w^{L1}_{n}}^{k^{\pi}_{n}} x dx = \frac{(k^{\pi}_{n})^{2}-(w^{L1}_{n})^{2}}{2}
> 0.
\end{displaymath}
Since among utilitarian equilibria the minimum value is minimised at the
libertarian equilibria when $N=2$, it follows by
induction that among utilitarian equilibria the minimum value is minimised at
the libertarian equilibria for $n\geq 2$.
By symmetry $v^{L1}_{n}+w^{L1}_{n} = v^{L2}_{n}+w^{L2}_{n}$.

We now show that the libertarian strategies are the only globally optimal utilitarian
strategies for $N\geq 3$. From the analysis of the two
stage game $v^{\pi}_{2}+w^{\pi}_{2}=\frac{3b^{2}}{2}-\frac{b}{2}+1$
for any utilitarian equilibrium.
Suppose $v^{L1}_{n}+w^{L1}_{n} > v^{\pi}_{n}+w^{\pi}_{n}$.
From the conditions for a utilitarian equilibrium $\pi$, it
follows that
\begin{eqnarray*}
v^{\pi}_{n+1}+w^{\pi}_{n+1} & = &
 \int_{0}^{k^{\pi}_{n}} (v^{\pi}_{n} + w^{\pi}_{n}) dx
+ \int_{k^{\pi}_{n}}^{1} (x+u_{n}) dx \\
v^{L1}_{n+1}+w^{L1}_{n+1} - (v^{\pi}_{n+1}+w^{\pi}_{n+1}) & = &
\int_{0}^{w^{L1}_{n}} [v^{L1}_{n}+w^{L1}_{n} - (v^{\pi}_{n}+w^{\pi}_{n})] dx
+\int_{w^{L1}_{n}}^{k^{\pi}_{n}} [x+u_{n} - (v^{\pi}_{n}+w^{\pi}_{n})] dx \\
& > & \int_{w^{L1}_{n}}^{k^{\pi}_{n}} [x+u_{n} - (v^{\pi}_{n}+w^{\pi}_{n})] dx \\
& = & \int_{w^{L1}_{n}}^{k^{\pi}_{n}} [x+u_{n} - (v^{L1}_{n}+w^{L1}_{n})
+ v^{L1}_{n}+w^{L1}_{n} - (v^{\pi}_{n}+w^{\pi}_{n})  ] dx > 0.\\
\end{eqnarray*}
This inequality follows from the induction assumption
$v^{L1}_{n}+w^{L1}_{n} - (v^{\pi}_{n}+w^{\pi}_{n})>0$, together with
 $v^{L1}_{n}<u_{n}$. It can be shown that
$v^{L1}_{3}+w^{L1}_{3} - (v^{\pi}_{3}+w^{\pi}_{3})>0$ using a similar argument for
$n=2$ (the
first inequality in the argument becomes an equality).
It follows by induction that for $N\geq 3$ the
libertarian equilibria are the only utilitarian equilibria which are globally optimal in the
sense of the utilitarian criterion.

\section{Final remarks}
In his recent paper, Garnaev \cite{gar04:staff} has extended the game model introduced in
 Baston and Garnaev \cite{basgar04:staff} to consider the situation where three
 skills of the candidate are taken into account. The proposed solutions to Garnaev's problem are
 Nash equilibria and Stackelberg strategies, as in \cite{basgar04:staff}, and
 these solutions are derived in his paper. One can also construct correlated equilibria for this model, which
 will be the subject of further investigation.

\bibliographystyle{plain}

\end{document}